\documentclass[12pt]{article}

\usepackage[utf8]{inputenc}
\usepackage[T1]{fontenc}
\usepackage[english]{babel}

\usepackage[a4paper,vmargin={2.5cm,2.5cm},hmargin={2cm,2cm}]{geometry}
\usepackage{mathpazo}
\usepackage[font=sf, labelfont={sf,bf}, margin=1cm]{caption}

\usepackage[pdftex]{hyperref}
\usepackage[expansion]{microtype}
\usepackage[pdftex]{color,graphicx}

\usepackage{amsmath}

\usepackage{amsfonts,amssymb,amsthm,mathrsfs}
\usepackage{stackrel}

\theoremstyle{plain}
\newtheorem*{theorem}{Theorem}

\theoremstyle{remark}
\newtheorem*{remark}{Remark}

\newcommand{\CC}{\mathrm{Cat}}

\newcommand{\HH}{\mathcal{H}}

\def\presuper#1#2%
  {\mathop{}%
   \mathopen{\vphantom{#2}}^{#1}%
   \kern-\scriptspace%
   #2}

\hypersetup{
    pdfauthor   = {N. Enriquez, L. Ménard},%
    pdftitle    = {Asymptotic expansion of the expected spectral measure of Wigner matrices}}

\begin{document}

\title{\bf Asymptotic expansion of the expected spectral measure of Wigner matrices}
\author{Nathanaël Enriquez, Laurent Ménard
}
\date{}
\maketitle


\begin{abstract}
We compute an asymptotic expansion with precision $1/n$ of the moments of the expected empirical spectral measure of Wigner matrices of size $n$ with independent centered entries. We interpret this expansion as the moments of the addition of the semi-circle law and $1/n$ times an explicit signed measured with null total mass. This signed measure depends only on the second and fourth moments of the entries.
\end{abstract}

\noindent{\bf {\textsc MSC 2010 Classification}:} 60B20.\\
\noindent{\bf Keywords:} Random matrices; moments method.

\section{Introduction and main result}

Let $n \geq 1$ be an integer. In this note, we will be interested in the so called Wigner ensemble of $n \times n$ matrices of the form
\begin{equation}
\label{eq:Wigner}
X^{(n)} = \left( X_{i,j} ^{(n)} \right)_{1 \leq i,j \leq n}, \quad \quad X_{i,j} ^{(n)} = \frac{1}{\sqrt{n} \sigma^2 } W_{i,j}
\end{equation}
where $\left( W_{i,j} \right)_{i \geq j \geq 1}$ are independent random variables with mean $0$ and moments of all order. We will consider two cases depending on whether the random variables are real or complex:
\begin{itemize}
\item {\bf Real case:} we set $r = 1$ and the random variables are such that $W_{j,i} = W_{i,j}$ and
\begin{equation}
\label{eq:real}
E[ W_{i,j} ] = 0, \quad E[ W_{i,j}^2] = 
\begin{cases}
\sigma^2 & \text{if $i \neq j$}\\
s^2 & \text{if $i=j$}
\end{cases},
\quad E[ W_{i,j} ^4 ] = \alpha \text{ if $i \neq j$}.
\end{equation}

\item {\bf Complex case:} we set $r = 0$ and the random variables are such that $W_{j,i} = \overline{W_{i,j}}$ and
\begin{equation}
\label{eq:complex}
E[ W_{i,j} ] = E[ W_{i,j}^2 ] = 0, \quad E[ \left| W_{i,j} \right|^2] = 
\begin{cases}
\sigma^2 & \text{ if $i \neq j$}\\
s^2 & \text{if $i=j$}
\end{cases},
\quad E[ \left| W_{i,j} \right| ^4 ] = \alpha \text{ if $i \neq j$}.
\end{equation}
\end{itemize}

Among the real and complex cases, one can find the Gaussian Orthogonal Ensemble (GOE, corresponding to $r = 1$, $s^2 = 2$, $\sigma^2 = 1$ and $\alpha = 3$) and the Gaussian Unitary Ensemble (GUE, corresponding to $r = 0$, $s^2 = 1$, $\sigma^2 = 1$ and $\alpha = 2$).

\bigskip

Our main purpose is to study the spectral measure
\begin{equation}
\label{eq:spectralmeasure}
\mu_n = \frac{1}{n} \sum_{\lambda \in \mathrm{Sp} (X^{(n)})} \delta_{\lambda}.
\end{equation}
of the matrices $X^{(n)}$.
More precisely, we compute an asymptotic expansion in $n$ of the expected moments
\begin{equation}
\label{eq:momentd}
m_k(n) = E \left[ \int x^k d \mu_n(x) \right]
\end{equation}
of the spectral measure $\mu_n$:

\begin{theorem}
For every $k \geq 0$, as $n \to \infty$, one has
\begin{equation}
\label{eq:mainresult}
m_k(n) = \int x^k d \mathrm{sc}(x) + \frac{1}{n} \int x^k d\nu(x) + o\left(\frac{1}{n} \right)
\end{equation}
where $\mathrm{sc}$ denotes the semi-circle measure having density
\[
\frac{1}{2 \pi} \sqrt{4 -x^2} \mathbf{1}_{ \{x \in [-2,2]\}}, 
\]
and $\nu$ is the signed measure with $0$ total mass given by:
\begin{align*}
& \frac{r}{2}  \left( \frac{1}{2} \left( \delta_{\{ x =2 \}} + \delta_{ \{ x = -2\}} \right) -\frac{\mathbf{1}_{\{ x \in [-2,2] \}} dx}{\pi \sqrt{4-x^2}}\right)\\
& + \frac{1}{2} \left(
\left( \frac{\alpha}{\sigma^4} - (2+ r) \right) x^4
+ \left( \frac{s^2}{\sigma^2} - 4 \frac{\alpha}{\sigma^4} +7 + 3 r\right) x^2
+ 2 \left( \frac{\alpha}{\sigma^4} - \frac{s^2}{\sigma^2} -1 \right)\right)
\frac{\mathbf{1}_{\{ x \in [-2,2] \}} dx}{\pi \sqrt{4-x^2}}.
\end{align*}
\end{theorem}

\begin{remark}
For the GUE, the measure $\nu$ is null, and for the GOE it is simply
\[
\frac{1}{2}  \left( \frac{1}{2} \left( \delta_{\{ x =2 \}} + \delta_{\{ x = -2 \}} \right) -\frac{\mathbf{1}_{\{ x \in [-2,2] \}} dx}{\pi \sqrt{4-x^2}}\right).
\]
\end{remark}

\bigskip

Before going through the proof of the theorem, let us make some comments. Asymptotic expansions of the expected spectral measure are of crucial interest for central limit theorems. Indeed, they are usually stated as follows for Wigner matrices: if $f$ is smooth enough, the random quantity
\[
n \left( \int f(x) d \mu_n (x) - \mathbb{E} \int f(x) d \mu_n (x)  \right)
\]
converges to a Gaussian random variable as $n \to \infty$ as proved by Lytova and Pastur in \cite{LP}. Hence it is important to know the centering term with a precision of order $o(1/n)$.

For the special case of the GUE, the measure $\nu$ is null. This agrees with the full genus expansion of $m_k(n)$ computed by Brezin Itzykson Parisi and Zuber \cite{BIPZ} and Harer and Zagier \cite{HZ}. For the GOE, Ledoux \cite{L} computed an analog of this full expansion which is also in agreement with our result.
The full expansion of $m_k(n)$ is not known for other ensembles. However Johansson \cite{J} computed the term of order $1/n$ for invariant matrix ensembles with a potential. This framework coincides with the case of matrices with independent entries only for the GUE and GOE.

The first order in the asymptotic expansion of the expected Stieltjes transform of $\mu_n$ is computed by Khorunzhy Khoruzhenko and Pastur \cite{KKP} for general matrices with independent entries. Their result is valid only at some finite distance from the real line. This analytic approach does not yield an asymptotic expansion for $m_k(n)$ and is radically different from the combinatorical method used in our computations of these moments.

\bigskip

Next section is devoted to the proof which enumerates some classes of closed paths which are fancy variants of the contour functions of trees. They are reminiscent of the combinatorics of trees and maps involved in the Harer Zagier expansion of $m_k(n)$ for the GUE.

\section{Proof of the theorem}

We start by writing the moments
\begin{align*}
m_k(n) & = \mathbb E \left[ \int x^k d \mu_n(x) \right]\\
& = \frac{1}{n^{1+\frac{k}{2}} \sigma^{k}} \sum_{i_1,\ldots, i_k = 1}^n \mathbb{E} \left[ W_{i_1,i_2} W_{i_2,i_3} \ldots W_{i_k,i_1}\right].
\end{align*}
As it is done classically, we regroup words $i_1 i_2 \ldots i_k i_1$ in equivalence classes (we say that two words are in the same class if they coincide up to a permutation of the letters of the alphabet $\{ 1, \ldots , n \}$). The main reason to do so is that two equivalent words $i_1 i_2 \ldots i_k i_1$ and $i'_1 i'_2 \ldots i'_k i'_1$ have the same contribution in $m_k(n)$:
\[
\mathbb{E} \left[ W_{i_1,i_2} W_{i_2,i_3} \ldots W_{i_k,i_1}\right] = \mathbb{E} \left[ W_{i'_1,i'_2} W_{i'_2,i'_3} \ldots W_{i'_k,i'_1}\right].
\]
If $c$ is an equivalence class of words, we denote by $\mathbb{E} \left[ W_c \right]$ the common contribution of its elements ({\it i.e.} the quantity $\mathbb{E} \left[ W_{i_1,i_2} W_{i_2,i_3} \ldots W_{i_k,i_1}\right]$ for one of its elements $i_1 i_2 \ldots i_k i_1$). We will often see a word $i_1 i_2 \ldots i_k i_1$ as a closed spanning walk with $k$ steps on a connected graph with vertex set $ \{ i_1, \ldots , i_k \}$ and edge set $\{ \{i_1,i_2\}, \{i_2,i_3 \} , \ldots , \{i_k, i_1\} \}$ (these graphs have no multiple edges but can have self-loops). Notice that the graphs and walks associated to two equivalent words are isomorphic.

\bigskip

Fix $k,v$ and $e$ some integers, we denote by $\mathscr{C}(k,v,e)$ the set of all equivalence classes of words $i_1 i_2 \ldots i_k i_1$ on the alphabet $\{ 1, \ldots , n \}$ such that
\[
\sharp \{ i_1, \ldots , i_k \} = v \quad \text{and} \quad \sharp \{ \{i_1,i_2\}, \{i_2,i_3 \} , \ldots , \{i_k, i_1\} \} = e.
\]
If $\mathscr{C}(k,v,e) \neq \emptyset$ satisfies $\mathbb E [W_c] \neq 0$ for some $c \in \mathscr{C}(k,v,e)$, one always has
\[
v \leq e+1 \quad \text{and} \quad e \leq \lfloor k/2 \rfloor.
\]
The first inequality is due to the fact that any connected graph with $v$ vertices has at least $v-1$ edges. The second inequality comes from the fact that the random variables $W_{i,j}$ are centered and therefore every edge of the graph defined by $c$ has to appear at least twice in a word belonging to $c$.

Now notice that every $c \in \mathscr{C}(k,v,e)$ is the equivalence class of exactly $n(n-1) \cdots (n-v+1)$ different words. This gives the expression
\begin{equation}
\label{eq:mk}
m_k(n) = \frac{1}{n^{1+\frac{k}{2}} \sigma^k} \sum_{e = 1}^{\lfloor k/2 \rfloor} \sum_{v = 1}^{e+1} n (n-1) \cdots (n-v+1) \sum_{c \in \mathscr{C}(k,v,e)} \mathbb{E} \left[ W_c \right].
\end{equation}
We are now ready to compute the asymptotic expansion of $m_k(n)$ announced in the theorem. First, we prove that the momemts to the moments of the semi-circle law, next we prove that the magnitude of the next term in the expansion is $1/n$ and finally we compute the term of order $1/n$ in the expansion. We identify these terms as the moments of the measure $\nu$ given in the theorem by computing their Stieltjes transform.

\subsection{Limits of the moments}

Non vanishing terms in \eqref{eq:mk} are such that $v = k/2$, and therefore $e = k/2$ and $k$ must be even. In that case every $c \in \mathscr{C} (k,k/2 + 1 , k/2)$ is such that $\mathbb E [W_c ] = (\sigma^2)^{k/2}$ and 
\[
\sharp \mathscr{C} (k,k/2 + 1 , k/2) = \mathrm{Cat} (k/2)
\]
yielding the moments of the semi-circle law for $\lim{n \to \infty} m_k(n)$. This is basically Wigner's original proof of convergence to the semi-circle law.

\subsection{The leading vanishing term is of order $1/n$}

Fix $a >0$. In order for \eqref{eq:mk} to have a term of order $1/ n^a$, there must be $1 \leq v \leq \lfloor \frac{k}{2} \rfloor + 1$ such that one has:
\[
\frac{n^{v}}{n^{\frac{k}{2} + 1}} = \frac{1}{n^a}
\]
and therefore
\[
a = 1 + \frac{k}{2} - v \geq \frac{k}{2} - \lfloor k/2 \rfloor \geq
\begin{cases}
0 & \text{if $k$ even;}\\
1/2 & \text{if $k$ odd.}
\end{cases}
\]
Now suppose $k=2l +1$. There is a term of order $1/\sqrt n$ in \eqref{eq:mk} if and only if $v = l+1$ and $e = l$. This means that the graph spanned by the word $i_1 i_2 \ldots i_k i_1$ is a tree with $l$ edges. In addition, the word is a closed walk of length $2l+1$ on this tree, which is impossible. Therefore, the leading vanishing term in \eqref{eq:mk} is at least of order $1/n$.

Notice that this also means that if $k$ is odd, the leading term in \eqref{eq:mk} is of order smaller than or equal to $n^{-3/2}$.

\subsection{The term of order $1/n$}

We suppose now that $k = 2l$ is even and compute the contributions of order $1/n$ in
\begin{equation}
\label{eq:mk2}
m_{2l}(n) = \frac{1}{n^{1+l} \sigma^l} \sum_{v,e} n (n-1) \cdots (n-v+1) \sum_{c \in \mathscr{C}(2l,v,e)} \mathbb E [W_c]
\end{equation}
where the sum on $v,e$ is taken over $v$ and $e$ satisfying
\[
1 \leq v \leq e+1 \quad \text{and} \quad 1 \leq e \leq k.
\]
The values of $v$ and $e$ yielding terms of order $1/n$ in \eqref{eq:mk2} are:
\begin{enumerate}
\item $v = e+1 = l +1$;
\item $v = e+1 = l$;
\item $v = e = l$.
\end{enumerate}

\subsubsection*{Case 1: $v = e+1 = l +1$.}
In this case, the contribution of order $1/n$ is the product of the coefficient of $n^{v-1}$ in the product $n(n-1) \cdots (n-v+1)$ and $Cat(l)$ as computed above for the order $0$. The total contribution of these terms is therefore
\[
m_{2l}^{(1)} = - \frac{1}{n} \frac{l(l+1)}{2} \mathrm{Cat}(l).
\]

\subsubsection*{Case 2: $v = e+1 = l$.}

In this case, the graph defined by a class $c \in \mathscr{C}(2l,v,e)$ is a tree with $l-1$ edges. In opposition to the previous case and the computation of order $0$, where each edge of the tree was visited exactly twice, here, exactly one edge is visited four times and all the other edges are visited twice. Note also that because the graph is a tree, if an edge $ij$ is visited in that order, the next visit of this edge must be done in the reverse order $ji$.

Consequently, we have to enumerate sequences of the form 
$$\mathcal{S} = i_1 \mathbf{S}_1 i \,  j \mathbf{S}_2  j  \, i \mathbf{S}_3 i \, j \mathbf{S}_4 j \, i \mathbf{S}_5 i_1 $$
such that the sequences $i_1 \mathbf{S}_1 i \mathbf{S}_5 i_1$, $j \mathbf{S}_2 j$, $i \mathbf{S}_3 i$ and $j \mathbf{S}_4 j$ are the contour functions of disjoint trees with respectively $p_1$, $p_2$, $p_3$ and $p_4$ edges satisfying $p_1 + p_2 + p_3 + p_4 = k -2$.

\begin{figure}[t!]
\begin{center}
\includegraphics[width=5cm]{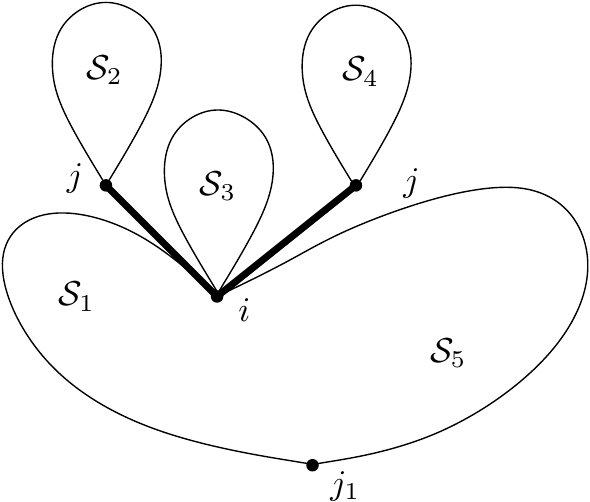}
\caption{\label{ERorder1_fig} Labeled tree with contour $\mathcal{S}$.}
\end{center}
\end{figure}

The sequence $i_1 \mathbf{S}_1 i \mathbf{S}_5 i_1$ corresponds to a rooted tree
with $p_1$ edges and a marked corner (adjacent to the vertex $i$) where the
trees corresponding to the three other sequences are inserted (see Figure
\ref{ERorder1_fig} for an illustration). Note that when $p_1 = 0$, the sequence
$\mathcal{S}$ boils down to $i \, j \mathbf{S}_2 j  \, i \mathbf{S}_3 i \, j
\mathbf{S}_4 j \, i$.

There are $(2p_1 +1) \CC(p_1)$ such trees with a marked corner. In that case every $c \in \mathscr{C} (2l, l , l-1)$ is such that $\mathbb E [W_c ] = (\sigma^2)^{l - 2} \alpha$ and 
\[
\sharp \mathscr{C} (2l, l, l-1) = \sum_{p_1 +\cdots + p_4 = l-2} (2p_1 +1) \CC(p_1) \cdots \CC(p_4).
\]
The contribution to the term of order $1/n$ in the asymptotic expansion of the moment $m_{2l}$ of this case is then:
\[
m_{2l}^{(2)} = \frac{1}{n} \frac{\alpha}{\sigma^4} \sum_{p_1 +\cdots + p_4 = l-2} (2p_1 +1) \CC(p_1) \cdots \CC(p_4).
\]

\subsubsection*{Case 3: $v = e = l$.}

In contrast with previous cases, there are two kinds of equivalent classes in  $\mathscr{C} (2l, l , l)$ leading to two different contributions $\mathbb E [ W_c]$. We detail each kind in the following.

\subsubsection*{First kind: graphs with a self-loop.}

Suppose that for some $1 \leq j \leq 2l$ we have $i_j = i_{j+1}$. We have to enumerate sequences of the form
\[
i_1 \mathbf{S}_1 i \,  i \mathbf{S}_2  i  \, i \mathbf{S}_3  i_1
\]
and such that the sequences $i_1 \mathbf{S}_1 i \mathbf{S}_3 i_1$ and $i \mathbf{S}_2 i$ are the contour functions of trees with respectively $p_1$ and $p_2$ edges satisfying $p_1 + p_2 = l - 1$.

For this kind of equivalence class, $\mathbb E [W_c ] = (\sigma^2)^{l - 1} s^2$ and there are
\[
\sum_{p_1 + p_2 = l-1} (2p_1 +1) \CC(p_1) \CC(p_2)
\]
such equivalent classes.
Hence, the contribution to the term of order $1/n$ in the asymptotic expansion of the moment $m_{2l}$ of this case is:
\[
m_{2l}^{(3)} = \frac{1}{n} \frac{s^2}{\sigma^2} \sum_{p_1 + p_2 = l-1} (2p_1 +1) \CC(p_1) \CC(p_2).
\]

\subsubsection*{Second kind: graphs with no self-loop.}

In that case, the underlying graph spanned by the walk $i_1i_2 \ldots i_{2l}i_1$ is a tree with $l-1$ edges and an additional edge between two of its vertices forming a cycle. The graph thus consits on a cycle of length $p \geq 3$ and trees attached to vertices of the cycle. We denote by $c_1,\ldots,c_p$ the vertices of the cycle.

First, notice that since $e=l$, each edge is visited exactly twice by the walk. We distinguish two cases:
\begin{enumerate}
\item edges of the cycle are visited one way (twice from $i$ to $j$);
\item edges of the cycle are visited both ways (once from $i$ to $j$ and once from $j$ to $i$ as in the previous cases).
\end{enumerate}

\begin{figure}[t!]
\begin{center}
\includegraphics[width=15cm]{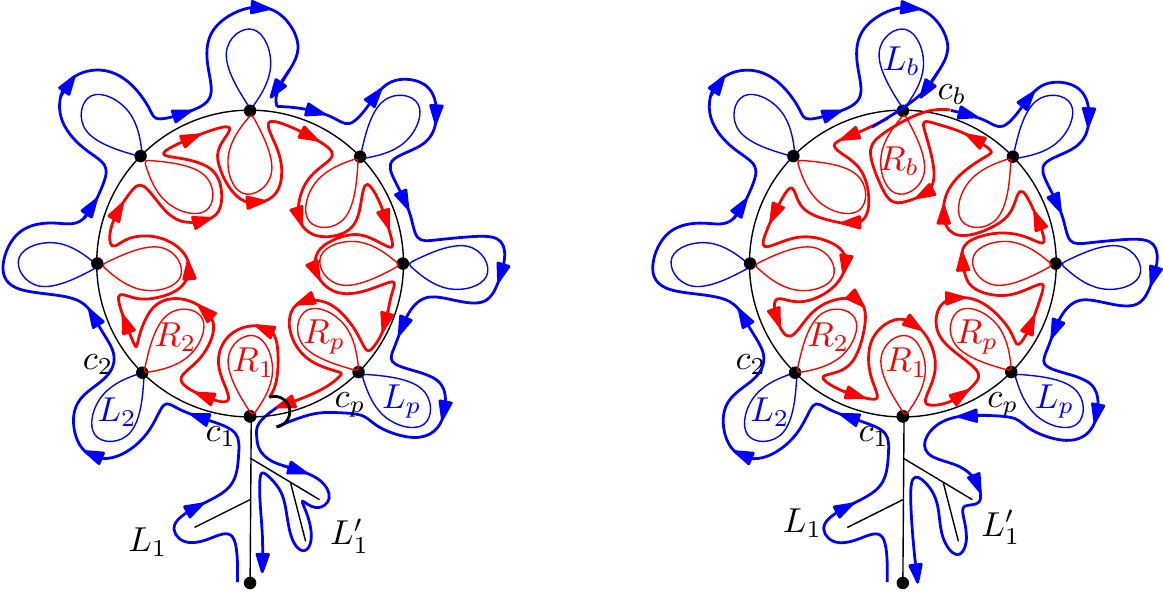}
\caption{\label{fig:cycle} A graph with one cycle. Left: edges of the cycle are visited one way only. Right: edges of the cycle are visited both ways.}
\end{center}
\end{figure}

The first case contributes only for real matrices and corresponds to sequences of the form
\[
i_1 L_1 c_1 c_2 L_2 c_2 c_3 \ldots c_p L_p c_p c_1 R_1 c_1 c_2 R_2 c_2 c_3 \ldots c_p R_p c_p c_1 L'_1 i_1
\]
such that the sequences $i_1 L_1 c_1 L'_1 i_1$, $(c_jL_jc_{j+1})_{2 \leq j \leq p}$ and $(c_jR_jc_{j+1})_{1 \leq j \leq p}$ (with $c_{p+1} = c_1$) code disjoint trees, see Figure \ref{fig:cycle} (left) for an illustration.

For this type of equivalence class, $\mathbb E [W_c ] = r (\sigma^2)^{l}$ and there are
\[
\sum_{p = 3}^l\sum_{\substack{l_1 + \cdots + l_p\\ +r_1 + \cdots r_p = l - p}} (2 l_1 +1) \prod_{i=1}^p \CC(l_i) \CC(r_i)
\]
such equivalent classes.
Hence, the contribution to the term of order $1/n$ in the asymptotic expansion of the moment $m_{2l}$ of this subcase is:
\[
\frac{r}{n} \sum_{p = 3}^l\sum_{\substack{l_1 + \cdots + l_p\\ +r_1 + \cdots r_p = l - p}} (2 l_1 +1) \prod_{i=1}^p \CC(l_i) \CC(r_i).
\]

\bigskip

Finally, let us consider the second case, where edges of the cycle are visited both ways. This corresponds to sequences of the form:
\begin{align*}
i_1 L_1 c_1 c_2 L_2 c_2 c_3 & \ldots L_{b-1}c_{b-1} c_b L_b c_b c_{b-1} R_{b-1} c_{b-1} \ldots c_2 c_1 R_1 c_1\\
& c_p L_p c_p c_{p-1} \ldots L_{b+1} c_{b+1} c_{b} R_{b} c_{b} c_{b+1} R_{b+1} \ldots c_p R_p c_p c_1 L'_1 i_1
\end{align*}
such that the sequences $i_1 L_1 c_1 L'_1 i_1$, $(c_jL_jc_{j+1})_{2 \leq j \leq p}$ and $(c_jR_jc_{j+1})_{1 \leq j \leq p}$ (with $c_{p+1} = c_1$) code disjoint trees, see Figure \ref{fig:cycle} (right).

For this type of equivalence class, $\mathbb E [W_c ] = (\sigma^2)^{l}$ and there are
\[
\sum_{p = 3}^l\sum_{\substack{l_1 + \cdots + l_p\\ +r_1 + \cdots r_p = l - p}} p (2 l_1 +1) \prod_{i=1}^p \CC(l_i) \CC(r_i)
\]
such equivalent classes.
Hence, the contribution to the term of order $1/n$ in the asymptotic expansion of the moment $m_{2l}$ of this subcase is:
\[
\frac{1}{n} \sum_{p = 3}^l\sum_{\substack{l_1 + \cdots + l_p\\ +r_1 + \cdots r_p = l - p}} p (2 l_1 +1) \prod_{i=1}^p \CC(l_i) \CC(r_i).
\]

\bigskip

Therefore, the total contribution to the term of order $1/n$ in the asymptotic expansion of the moment $m_{2l}$ of this case is:
\[
m_{2l}^{(4)} = \frac{1}{n} \sum_{p = 3}^l (p+r) \sum_{\substack{l_1 + \cdots + l_p\\ +r_1 + \cdots r_p = l - p}} (2 l_1 +1) \prod_{i=1}^p \CC(l_i) \CC(r_i).
\]

\subsection{Generating series}

We are now interested in the computation of the generating series of the terms of order $1/n$ where the generating series of the Catalan numbers will play a central role. We denote this series by
\[
T(x) = \sum_{k \geq 0} \CC(k) x^k.
\]
It has a radius of convergence $1/4$ and satisfies the following equalities:
\[
T(x) = 1 + x (T(x))^2 = \frac{1}{1 - xT(x)} = \frac{1 - \sqrt{1 - 4x}}{2x}.
\]
Our computations will also involve the first two derivatives of $T$, and
some basic computations yield the following useful equations:
\begin{align*}
T' & = \frac{T^3}{1 - xT^2};\\
T'' &= \frac{2T^5}{(1 - xT^2)^2} + \frac{2T^5}{(1 - xT^2)^3}.
\end{align*}

We start with the generating series of $m_{2k}^{(1)}$:
\begin{align*}
S_1(x) & = \sum_{k \geq 0} m_{2k}^{(1)} x^k = - \frac{1}{2} \sum_{k \geq 0} k(k+1) \CC(k) x^k\\
& = - \frac{1}{2} x^2 T''(x) - x T'(x)\\
& = - \frac{xT^3}{(1-xT^2)^3}.
\end{align*}

Similarly for $m_{2k}^{(2)}$:
\begin{align*}
S_2(x) & = \sum_{k \geq 0} m_{2k}^{(2)} x^k =  \frac{\alpha}{\sigma^4} \sum_{k \geq 2}
\sum_{p_1 +\cdots + p_4 = k-2} (2p_1 +1) \CC(p_1) \cdots \CC(p_4) x^k\\
& = \frac{\alpha}{\sigma^4} x^2 \left( \sum_{p_1 \geq 0} (2p_1 + 1) \CC(p_1) x^{p_1}\right) \left( \sum_{p_2 \geq 0} \CC(p_2) x^{p_2}\right)^3\\
& = \frac{\alpha}{\sigma^4} x^2 (2xT' +T) T^3\\
& = \frac{\alpha}{\sigma^4} \frac{x^2 T^5}{1-xT^2}.
\end{align*}

And for $m_{2k}^{(3)}$:
\begin{align*}
S_3(x) & = \sum_{k \geq 0} m_{2k}^{(3)} x^k =  \frac{s^2}{\sigma^2} \sum_{k \geq 1}
\sum_{p_1 + p_2 = k-1} (2p_1 +1) \CC(p_1) \CC(p_2) x^k\\
& = \frac{s^2}{\sigma^2} x \left( \sum_{p_1 \geq 0} (2p_1 + 1) \CC(p_1) x^{p_1}\right) \left( \sum_{p_2 \geq 0} \CC(p_2) x^{p_2}\right)\\
& = \frac{s^2}{\sigma^2} x (2xT' +T) T\\
& = \frac{s^2}{\sigma^2} \frac{x T^3}{1-xT^2}.
\end{align*}

And finally for $m_{2k}^{(4)}$:
\begin{align*}
S_4(x) & = \sum_{k \geq 0} m_{2k}^{(4)} x^k =  \sum_{k \geq 3}
\sum_{p = 3}^k (p+r) \sum_{\substack{l_1 + \cdots + l_p\\ +r_1 + \cdots r_p = k - p}} (2 l_1 +1) \prod_{i=1}^p \CC(l_i) \CC(r_i) x^k\\
& = \sum_{p \geq 3} (p+r) x^p T^{2p -1} (2xT'+T) = \frac{T}{1-xT^2} \sum_{p \geq 3} (p+r) x^p T^{2p}\\
& = \frac{x^3 T^7}{(1-xT^2)^3} + (2 +r) \frac{x^3 T^7}{(1-xT^2)^2}.
\end{align*}

The total generating function of terms of order $1/n$ is thus
\begin{align*}
S(x) & = S_1(x) + S_2(x) + S_3(x) + S_4(x)\\
& = \frac{x T^3}{(1-xT^2)^3} (x^2T^4 -1) + (2 +r) \frac{x^3 T^7}{(1-xT^2)^2} + \frac{x T^3}{1-xT^2} \left( \frac{\alpha}{\sigma^4} x T^2 + \frac{s^2}{\sigma^2} \right)\\
& =  - \frac{ x T^4}{(1-xT^2)^2} + (2 +r) \frac{x^3 T^7}{(1-xT^2)^2} + \frac{x T^3}{1-xT^2} \left( \frac{\alpha}{\sigma^4} x T^2 + \frac{s^2}{\sigma^2} \right)\\
& =  - \frac{ x T^4}{(1-xT^2)^2} +2 \frac{x^3 T^7}{(1-xT^2)^2} + 2 \frac{x^2 T^5}{1-xT^2} + \frac{x T^3}{1-xT^2} \\
& \quad \quad + r \frac{x^3 T^7}{(1-xT^2)^2} + \frac{x T^3}{1-xT^2} \left( \left(\frac{\alpha}{\sigma^4} -2 \right) xT^2 + \frac{s^2}{\sigma^2} -1 \right).
\end{align*}
We gathered terms in the last equality so that the last four terms vanish in the GUE setting for which $r = 0$, $\alpha / \sigma^4 = 2$ and $s^2 / \sigma^2 =1$. The first four terms therefore correspond to the generating series of the terms of order $1/n$ for the GUE which we know to be null.
It is easy to verify that indeed
\[
- \frac{ x T^4}{(1-xT^2)^2} +2 \frac{x^3 T^7}{(1-xT^2)^2} + 2 \frac{x^2 T^5}{1-xT^2} + \frac{x T^3}{1-xT^2} = 0
\]
and therefore
\begin{align}
\label{eq:genseries}
S(x) & = r \frac{x^3 T^7}{(1-xT^2)^2} + \frac{x T^3}{1-xT^2} \left( \left(\frac{\alpha}{\sigma^4} -2 \right) xT^2 + \frac{s^2}{\sigma^2} -1 \right) \notag\\
& = r \left( \frac{x^3 T^7}{(1-xT^2)^2} + \frac{x T^3}{1-xT^2} (xT^2 + 1) \right) + \frac{x T^3}{1-xT^2} \left( \left(\frac{\alpha}{\sigma^4} -2 -r \right) xT^2 + \frac{s^2}{\sigma^2} -1 -r \right) \notag \\
& = r \frac{x T^3}{(1-xT^2)^2}\left( x^2 T^4 + (xT^2 + 1) (1-xT^2) \right) + \frac{x T^3}{1-xT^2} \left( \left(\frac{\alpha}{\sigma^4} -2 -r \right) xT^2 + \frac{s^2}{\sigma^2} -1 -r \right).
\end{align}

\bigskip

Let
\[
\mathcal{H} (z) = \frac{1}{z} T \left( \frac{1}{z^2} \right) = \frac{1}{2} \left( z - \sqrt{z^2 -4} \right)
\]
be the Stieltjes transform of the semi-circle law. It satisfies:
\[
1 - \mathcal{H}^2 = \HH \sqrt{z^2 -4}.
\]

A simple computation using \eqref{eq:genseries} yields
\begin{align*}
\tilde{ \mathcal{H} } (z) &= \frac{1}{z} S \left( \frac{1}{z^2} \right) \\
&= r \frac{\mathcal H}{z^2 - 4} + \frac{\HH^2}{\sqrt{z^2 - 4}} \left( \left( \frac{\alpha}{\sigma^4} - 2 -r \right) \HH^2 + \frac{s^2}{\sigma^2} - 1 - r \right)\\
& = \frac{r}{2} \left( \frac{z}{z^2 - 4} - \frac{1}{\sqrt{z^2 - 4}} \right) + \frac{\HH^2}{\sqrt{z^2 - 4}} \left( \left( \frac{\alpha}{\sigma^4} - 2 -r \right) \HH^2 + \frac{s^2}{\sigma^2} - 1 - r \right)\\
&= \frac{r}{2}
\left( \frac{1}{2} \left( \frac{1}{z-2} + \frac{1}{z+2} \right) - \frac{1}{\sqrt{z^2 - 4}}\right)
+ \frac{\HH^2}{\sqrt{z^2 - 4}} \left( \left( \frac{\alpha}{\sigma^4} - 2 -r\right) \HH^2 + \frac{s^2}{\sigma^2} - 1 -r\right)
\end{align*}
which is the Stieltjes transform of the measure $\nu$ in the Theorem.

\bigskip

\noindent {\bf Acknowledgements:} N.E is partially supported by ANR PRESAGE (ANR-11-BS02-003) and L.M. is partially supported by ANR GRAAL (ANR-14-CE25-0014).

\addcontentsline{toc}{section}{References}
\bibliographystyle{abbrv}
\bibliography{matrices}

\noindent \textsc{Nathana\"el Enriquez} \verb|nenriquez@u-paris10.fr|, \\
\textsc{Laurent M\'enard} \verb|laurent.menard@normalesup.org|\\
Universit\'e Paris-Ouest Nanterre\\
Laboratoire Modal'X\\ 200 avenue de la R\'epublique\\
92000 Nanterre (France).

\end{document}